\documentclass[10pt,english]{smfart}
\usepackage[T1]{fontenc}
\IfFileExists{french.ldf}{\usepackage[french,english]{babel}}{%
  \usepackage[english]{babel}\babelprovide[import]{french}}
\usepackage{latexsym,amscd,color}
\usepackage{amsmath,amsfonts,amssymb,mathrsfs}
\usepackage{enumerate,euscript}
\usepackage{amssymb,url,xspace,smfthm}
\IfFileExists{smfhyperref.sty}{\usepackage{hyperref}}{}
\input xy
\xyoption{all}

\newcommand{\BibTeX}{{\scshape Bib}\kern-.08em\TeX}

\newcommand{\T}{\S\kern .15em\relax }
\newcommand{\AMS}{$\mathcal{A}$\kern-.1667em\lower.5ex\hbox
        {$\mathcal{M}$}\kern-.125em$\mathcal{S}$}

\DeclareMathOperator{\spm}{Spm}

\DeclareMathOperator{\spec}{Spec}

\renewcommand{\P}{\mathbb{P}}
\newcommand{\wmu}{\widehat{\mu}}

\newcommand{\Q}{\mathbb{Q}}

\newcommand{\adeg}{\widehat{\deg}}

\newcommand{\p}{\mathfrak{p}}

\newcommand{\sE}{\mathcal{E}}

\renewcommand{\O}{\mathcal{O}}
\renewcommand{\L}{\overline{L}}
\newcommand{\sm}{\mathfrak{m}}

\newcommand{\ndot}{\raisebox{.4ex}{.}}
\newcommand{\PP}{\mathbb P}

\newcommand{\FF}{\mathbb F}
\newcommand{\OO}{\mathcal O}
\newcommand{\fp}{\mathfrak p}
\newcommand{\cX}{\mathscr X}

\newcommand{\Np}{N(\fp)}

\DeclareMathOperator{\rk}{rk}
\DeclareMathOperator{\length}{length}
\DeclareMathOperator{\coker}{coker}

\title[\texorpdfstring{Corrigendum to ``Explicit uniform estimation of rational points II''}{Corrigendum to Explicit uniform estimation of rational points II}]{Corrigendum to ``Explicit uniform estimation of rational points II. Hypersurface coverings''}
\alttitle{Corrigendum à « Explicit uniform estimation of rational points II. Hypersurface coverings »}

\date{}
\author{Chunhui Liu}
\address{Institute for Advanced Study in Mathematics\\
Harbin Institute of Technology\\150001 Harbin\\P. R. China}
\email{chunhui.liu@hit.edu.cn}

\begin{document}
\def\smfbyname{}
\begin{abstract}
We correct an error in H.~Chen's treatment of non-reduced special fibers in his proof of a uniform bound for rational points of bounded height on integral plane curves over number fields. The correction completes the argument within its original Arakelov-geometric framework, using a multiplicity estimate together with slope inequalities. We obtain a corrected covering by auxiliary curves and recover the claimed counting bound. The constants are explicit and uniform over number fields of fixed degree.
\end{abstract}

\maketitle

\tableofcontents
\section{Introduction and preliminaries}
In this article, we correct an error in the proof of \cite[Theorem~5.1]{Chen2}, from which H. Chen deduced a bound for rational points on plane curves in \cite[Corollary~5.2]{Chen2}, in response to a question of Heath-Brown \cite[Question~27]{ARCC02}.
\subsection{Background}
Counting rational points of bounded height is a fundamental problem in arithmetic geometry. For a number field $K$ and an integral plane curve $X\subset\mathbb P^2_K$ of degree $\delta$, we use $S(X;B)$ and the relative height defined in \S\ref{fundamental setting}. In the classical setting over $\Q$, Heath-Brown asked for improved uniform bounds when the height bound equals the degree; see \cite[Question~27]{ARCC02}. A stronger form of this counting problem is the estimate
\[
 \#S(X;B)\ll_{K,\varepsilon}\delta^2B^{2/\delta+\varepsilon},\qquad B\geqslant2,\quad\varepsilon>0,
\]
uniformly over integral plane curves $X$. The present article concerns its weaker consequence \begin{equation}\label{intro:diagonal}\#S(X;\delta)\ll_{K,\varepsilon}\delta^{2+\varepsilon}.\end{equation}

In \cite{Chen1,Chen2}, H. Chen developed an Arakelov-geometric formulation of the determinant method, in which slope inequalities applied to evaluation maps produce auxiliary hypersurfaces. He claimed the estimate \eqref{intro:diagonal} in \cite[Corollary~5.2]{Chen2} as a consequence of \cite[Theorem~5.1]{Chen2}. However, the proof of that theorem applies the multiplicity estimate \cite[(30)]{Chen2} to special fibers which need not be reduced. In fact, that estimate is valid for reduced plane curves, but need not hold for non-reduced fibers. For example, for a rational prime $p$ the smooth conic $x_0^2+p x_1x_2=0$ over $\Q$ reduces modulo $p$ to a double line. Each of its $p+1$ rational points over $\mathbb F_p$ has multiplicity two, so the sum in \cite[(30)]{Chen2} is $2(p+1)$, whereas the asserted upper bound is $2$.

Binyamini, Cluckers and F.~Kato subsequently proved a sharper counting bound over all global fields in \cite[Theorem~1]{BCK2025}. For each fixed number field $K$ and the relative height used here, their result implies
\[
 \#S(X;B)\ll_K\delta^2B^{2/\delta}(\log(2B))^{12},\ B\geqslant2.
\]
Their theorem uses the absolute projective height; the displayed estimate follows by passing to relative heights, comparing the supremum and Fubini--Study norms, and enlarging the implied constant for bounded $B$. Since $(\log(2B))^{12}\ll_\varepsilon B^\varepsilon$, this gives the stated bound for general $B\geqslant2$.

\subsubsection{}
We prove the following corrected covering estimate over a number field $K$. The explicit constant in the upper bound is modified and depends only on $[K:\Q]$.
\begin{theo}\label{th:main}
Let $X$ be an integral plane curve in $\mathbb P^2_K$ of degree $\delta\geqslant2$, and let $e<B<e^{\delta^2}$. Put $b_K=256[K:\Q]([K:\Q]+1)+2[K:\Q]+16$, which depends only on $[K:\Q]$. Then $S(X;B)$ is covered by at most
\begin{equation}\label{eq:bound}
 \exp(b_K\sqrt{\log B})+(\log B)^{\sqrt{\log B}+1}
\end{equation}
plane curves over $K$ of degree at most $\delta-1$. None contains $X$. In particular,
\[
 \#S(X;B)\leqslant\delta(\delta-1)
          \bigl(\exp(b_K\sqrt{\log B})+(\log B)^{\sqrt{\log B}+1}\bigr).
\]
\end{theo}
The following consequence is proved in Corollary~\ref{co:diagonal}. 
\begin{coro}
Let $X$ be an integral plane curve in $\mathbb P^2_K$ of degree $\delta\geqslant1$. Then, for every $\varepsilon>0$,
\[\#S(X;\delta)\ll_{[K:\Q],\varepsilon}\delta^{2+\varepsilon}.\]
The implied constant is explicit in terms of $[K:\Q]$ and $\varepsilon$; in particular, the bound is uniform over number fields of fixed degree.
\end{coro}
\subsubsection{}
We repair H.~Chen's proof within its original Arakelov-geometric framework by supplying the missing treatment of non-reduced fibers. The geometric separation of components of large multiplicity from a finite exceptional set controlled by intersection estimates also occurs in the work of Binyamini, Cluckers and F.~Kato; see \cite[proof of Proposition~3 and \S5.2]{BCK2025}.

In \cite{BCK2025}, the authors cover the exceptional points in each special fiber by a low-degree curve and take its union with the components of large multiplicity. We treat the two sets separately. Consider the rational points whose reductions have large multiplicity at every selected prime, and take $r$ as in \eqref{eq:parameters}. If a point reduces outside the selected components at $r$ or more primes, we fix $r$ such primes and its reduction at each. Lemma~\ref{lem:tube} gives one curve over $K$ containing all points with those prescribed reductions. Each remaining point reduces outside the selected components at most $r-1$ times. Proposition~\ref{prop:slope} covers all these points by one further curve.

Both arguments compare divisibility at prime ideals with height bounds. In \cite{BCK2025}, this comparison uses interpolation determinants formed in homogeneous coordinates. We formulate it through evaluation maps taking values in pullback Hermitian line bundles, using the slope inequalities and local determinant estimates of \cite{Chen1,Chen2}. These refinements complete H.~Chen's original argument without invoking the stronger counting theorem of \cite{BCK2025}.

A further feature is that the constants in Theorem~\ref{th:main} and Corollary~\ref{co:diagonal} are fully explicit
and depend on the number field only through $[K:\Q]$. The diagonal counting bound is therefore uniform over number fields of fixed degree, with an explicit dependence on $\varepsilon$. This uniformity is not part of the statement of \cite[Theorem~1]{BCK2025}, whose constant is allowed to depend on $K$ and is not given an explicit value.

\subsection{Fundamental setting}\label{fundamental setting}
Unless otherwise stated, we use the notation and conventions of \cite{Chen1,Chen2}. 
\subsubsection{}
Let $K$ be a number field and $\O_K$ its ring of integers. We denote the sets of finite and infinite places of $K$ by $M_{K,f}$ and $M_{K,\infty}$, respectively. We denote by $M_K$ the set of all places of $K$. Let $\spm\O_K$ be the set of maximal ideals of $\O_K$. For $\p\in\spm\O_K$, let $\O_{K,\p}$ be the localization of $\O_K$ at $\p$, let $\mathbb F_\p=\O_K/\p$ be its residue field, and put $N(\p)=\#\mathbb F_\p$. 

For $v\in M_K$, let $K_v$ and $\Q_v$ be the completions at $v$ and its restriction to $\Q$, respectively. We use the absolute value on $K_v$ extending the usual real or $p$-adic absolute value on $\Q_v$, with $|p|_p=p^{-1}$; at a complex place, we use the usual complex modulus.

For a Hermitian vector bundle $\overline E$ of positive rank over $\spec\O_K$, we denote by $\adeg(\overline E)$ its Arakelov degree and by
\[
 \wmu(\overline E)=\frac{\adeg(\overline E)}{[K:\Q]\rk E}
\]
its normalized slope, as in \cite[Notations~5--6]{Chen1}.

Let $\overline{\sE}=(\O_K^{\oplus 3},(\|\ndot\|_v)_{v\in M_{K,\infty}})$ be the Hermitian vector bundle over $\spec\O_K$, where for each $v\in M_{K,\infty}$, we define the norm
\[\begin{array}{rrcl}
\|\ndot\|_v:&K_v^3&\longrightarrow&\mathbb R\\
&(x_0,x_1,x_2)&\mapsto&\sqrt{|x_0|_v^2+|x_1|_v^2+|x_2|_v^2}.
\end{array}\]
\subsubsection{}
Let $\mathbb P(\sE)$ be the projective space over $\spec\O_K$ which parametrizes invertible quotients of $\sE$. We equip its universal quotient line bundle $\O(1)=\O_{\mathbb P(\sE)}(1)$ with the corresponding Fubini--Study metrics at the infinite places, and denote the resulting Hermitian line bundle by $\overline{\O(1)}$. Write $\sE_K=\sE\otimes_{\O_K}K$. Each point $P=[x_0:x_1:x_2]\in\mathbb P(\sE_K)(K)$ extends uniquely to a section $\mathcal P:\spec\O_K\to\mathbb P(\sE)$. For $\fp\in\spm\O_K$, write $P_\fp=\mathcal P(\fp)\in\mathbb P^2(\FF_\fp)$ for its reduction modulo $\fp$. We write
\begin{align*}
 h(P)&:=h_{\overline{\O(1)}}(\mathcal P)
       =\frac{\adeg(\mathcal P^*\overline{\O(1)})}{[K:\Q]}\\
 &=\sum_{v\in M_{K,f}}\frac{[K_v:\Q_v]}{[K:\Q]}
       \log\max_{0\leqslant i\leqslant2}\{|x_i|_v\}\\
 &\quad+\sum_{v\in M_{K,\infty}}\frac{[K_v:\Q_v]}{[K:\Q]}
       \log\sqrt{|x_0|_v^2+|x_1|_v^2+|x_2|_v^2}.
\end{align*}
Put $H(P)=\exp([K:\Q]h(P))$. This is the relative height of \cite[Notation~11]{Chen2}.

For a subscheme $X$ of $\mathbb P(\sE_K)$ and $B>0$, let $S(X;B)=\{P\in X(K)\mid H(P)\leqslant B\}$. 

For an integer $D\geqslant1$, let $E_D=H^0(\mathbb P(\sE),\O(D))$, where $\O(D)=\O(1)^{\otimes D}$. We equip $E_D$ at each infinite place with the John norm associated to the supremum norm on the corresponding complex projective plane, using the Fubini--Study metric, and denote the resulting Hermitian vector bundle by $\overline E_D$; see \cite[\S2.3]{Chen1}. We write \[r(2,D)=\rk (E_D)=\binom{D+2}{2}.\]

Let $X$ be an integral plane curve in $\mathbb P(\sE_K)$ of degree $\delta$, whose Zariski closure in $\mathbb P(\sE)$ is denoted by $\mathscr X$. Let $F_D$ be the saturation of the image of the restriction map
\[
 E_D\longrightarrow H^0(\cX,\O(D)|_{\cX})
\]
in its target. Following \cite[Notation~14]{Chen2}, let $\overline F_D$ be the Hermitian vector bundle on $\spec\O_K$ whose norms at the infinite places are the quotient norms induced by the John norms of $\overline E_D$. We write
\[
 r_1(D)=\rk (F_D)=\binom{D+2}{2}-\binom{D-\delta+2}{2},
\]
where the second binomial coefficient is taken to be zero when $D<\delta$; see \cite[\S3]{Liu2022d} for the rank formula.
\subsection*{Acknowledgement}
I would like to thank Prof. Huayi Chen for supervising my Ph.D. thesis, which was based on \cite{Chen1,Chen2}, and for many discussions on this subject. This work was completed during my visit to the Institut de Math\'ematiques de Jussieu--Paris Rive Gauche and Sorbonne Universit\'e. I am grateful to both institutions for their hospitality.

ChatGPT was used to assist in developing and checking parts of the mathematical argument and in revising the exposition.
\section{Counting multiplicities}
\markright{CORRIGENDUM TO H.~CHEN'S HYPERSURFACE COVERINGS}
In this section, we consider the number of points with large multiplicity on a plane curve. 
\subsection{Local multiplicities}
We recall the definition of multiplicity for local rings and for points on plane curves.
\subsubsection{}
Let $(A,\sm,\kappa)$ be a Noetherian local ring of Krull dimension $d$, whose maximal ideal is $\sm$ and residue field is $\kappa$. For $D\in\mathbb N$, put $H(D)=\dim_\kappa(\sm^{D}/\sm^{D+1})$. If $d\geqslant1$, we have the asymptotic expansion
\[H(D)=\frac{\mu}{(d-1)!}D^{d-1}+O(D^{d-2})\]
as $D\to\infty$, by \cite[Corollary 11.5, Theorem 11.14]{AM_commutative_algebra}. The positive integer $\mu$ is the \textit{multiplicity} of $A$. If $d=0$, we put $\mu=\length_A(A)$.

If $(A,\sm,\kappa)$ is a regular local ring, then $\mu=1$. 
\subsubsection{}
Let $X$ be a projective scheme over a field $k$. For $\xi\in X$, we denote by $\mu_\xi(X)$ the \textit{multiplicity} of the local ring $\O_{X,\xi}$. Its maximal ideal and residue field are denoted by $\sm_\xi$ and $\kappa(\xi)$, respectively.

We say that $\xi$ is \textit{regular} if $\O_{X,\xi}$ is a regular local ring, and \textit{singular} otherwise. For a plane curve locally defined by one equation, the multiplicity at a point is the order of that equation in the ambient regular local ring, so the point is regular if and only if its multiplicity is $1$. For a point $\xi$ of the ambient scheme outside $X$, we put $\mu_\xi(X)=0$. 
\subsection{Multiplicity estimates}
Let $\mathbb F_q$ be a finite field of cardinality $q$, and let $Y$ be an effective divisor of degree $\delta\geqslant1$ in $\mathbb P^2_{\mathbb F_q}$. Over the algebraic closure $\overline{\mathbb F}_q$ of $\mathbb F_q$, write $Y_{\overline{\mathbb F}_q}=\sum\limits_{i=1}^m a_i C_i$, where the $C_i$ are distinct integral plane curves and $a_i\geqslant1$. At geometric points, all multiplicities are taken after base change to $\overline{\mathbb F}_q$. At an $\mathbb F_q$-rational point, this agrees with the multiplicity before base change, since extending the coefficient field preserves the order of a local defining equation. 

\begin{lemm}\label{lem:weighted}
With the above notation, we have
\[
 \sum_{\xi\in Y(\overline{\mathbb F}_q)}
 \left(\mu_\xi(Y)^2-\sum_{i=1}^m a_i^2\mu_{\xi}(C_i)\right)
 \leqslant\deg(Y)^2-\sum_{i=1}^m a_i^2\deg(C_i).
\]
\end{lemm}
\begin{proof}
First, the summand vanishes at every point that is regular on a single component $C_i$ and lies on no other component. The remaining points belong to the singular loci of the $C_i$ or to their pairwise intersections. These sets are finite, since the $C_i$ are integral and distinct. The order of a product of local equations gives $\mu_\xi(Y)=\sum\limits_{i=1}^m a_i\mu_{\xi}(C_i)$, and degrees give $\deg(Y)=\sum\limits_{i=1}^m a_i\deg(C_i)$.

In the following two sums, the points are taken over $\overline{\mathbb F}_q$. For each $i=1,\ldots,m$, we have
\[
 \sum_{\xi\in C_i}\mu_{\xi}(C_i)(\mu_{\xi}(C_i)-1)\leqslant \deg(C_i)(\deg(C_i)-1)
\]
by \cite[Exercise 5--22]{Fulton2}. To apply this estimate in positive characteristic, some partial derivative of an irreducible equation for $C_i$ must be nonzero. Otherwise, that equation would be a $p$-th power over the perfect field $\overline{\mathbb F}_q$, where $p=\operatorname{char}(\mathbb F_q)$.

For $i\neq j$, the curves $C_i$ and $C_j$ have no common component. By B\'ezout's theorem and the local intersection multiplicity inequality \cite[Chap.~3, \S3, Property~(5), and Chap.~5, \S3]{Fulton2}, we obtain
\[\sum_{\xi\in C_i\cap C_j}\mu_{\xi}(C_i)\mu_{\xi}(C_j)\leqslant \deg(C_i)\deg(C_j).\]

Expanding the square and applying these two bounds yields
\begin{align}
& \sum_{\xi\in Y(\overline{\mathbb F}_q)}\left(\mu_\xi(Y)^2-\sum_{i=1}^m a_i^2\mu_{\xi}(C_i)\right)\notag\\
  \leqslant& \sum_{i=1}^m a_i^2\deg(C_i)(\deg(C_i)-1)+\sum_{i\ne j}a_i a_j \deg(C_i) \deg(C_j)\notag\\=&\deg(Y)^2-\sum_{i=1}^m a_i^2 \deg(C_i) ,\label{eq:weighted}
\end{align}
which proves the assertion. 
\end{proof}

\begin{prop}\label{lem:mult}
For an integer $w\geqslant1$, let $T$ be the reduced union of the $C_i$ such that $a_i\geqslant w/2$. Then $T$ is defined over $\mathbb F_q$, has degree at most $\lfloor 2\delta/w\rfloor$, and
\begin{equation}\label{eq:local}
 \#\left\{\xi\in Y(\overline{\mathbb F}_q)\smallsetminus T(\overline{\mathbb F}_q)\mid \mu_\xi(Y)\geqslant w\right\}
 \leqslant \left\lfloor\frac{2\delta^2}{w^2}\right\rfloor .
\end{equation}
The empty union is allowed and has degree zero.
\end{prop}

\begin{proof}
The $q$-power Frobenius automorphism $x\mapsto x^q$ of $\overline{\mathbb F}_q/\mathbb F_q$ permutes the components and preserves their multiplicities, so it permutes the selected components. In a fixed monomial order, normalize the product of their equations to have first nonzero coefficient one. The above Frobenius automorphism fixes this normalized equation, so its coefficients belong to $\mathbb F_q$. Its degree is \[\sum\limits_{a_i\geqslant w/2}\deg(C_i)\leqslant\frac{2}{w}\sum_{i=1}^m a_i\deg(C_i)=\frac{2\deg(Y)}{w}.\] The empty union is represented by the constant form $1$.

If $\xi\notin T$, every component $C_i$ through $\xi$ has $a_i<\frac{w}{2}$, so 
\[\sum_{i=1}^m a_i^2\mu_{\xi}(C_i)\leqslant\frac{w}{2}\mu_\xi(Y).\]
Each summand in \eqref{eq:weighted} is nonnegative, since it equals $\sum\limits_{i=1}^m a_i\mu_\xi(C_i)(\mu_\xi(Y)-a_i)$ and $\mu_\xi(Y)\geqslant a_i$ whenever $\xi\in C_i$. Lemma~\ref{lem:weighted} therefore bounds the sum over any subset of points by $\deg(Y)^2$.

At each point outside $T$ with $\mu_\xi(Y)\geqslant w$, we have
\[
 \mu_\xi(Y)^2-\sum_{i=1}^m a_i^2\mu_{\xi}(C_i)
 \geqslant\mu_\xi(Y)\left(\mu_\xi(Y)-\frac{w}{2}\right)
 \geqslant\frac{\mu_\xi(Y)^2}{2}\geqslant\frac{w^2}{2}.
\]
The number of these points is therefore at most $\displaystyle\frac{2\deg(Y)^2}{w^2}$, and we may take the floor because the cardinality is an integer. This proves \eqref{eq:local}. 
\end{proof}
\section{Evaluation maps and auxiliary curves}
In this section, we retain the notation of \S\ref{fundamental setting}.

For $D\geqslant1$ and a finite family of distinct points $(P_i)_{i\in I}$ in $\mathbb P(\sE_K)(K)$, let $\mathcal P_i$ be their extensions over $\O_K$. Their evaluation map is
\[
\begin{array}{rrcl}
\phi:&E_D&\longrightarrow&\bigoplus\limits_{i\in I}\mathcal P_i^*\O(D)\\
&s&\longmapsto&(\mathcal P_i^*s)_{i\in I}.
\end{array}
\]
The target carries the orthogonal direct sum of the pullback metrics. The kernel of $\phi_K=\phi\otimes_{\O_K}K$ consists of the degree $D$ forms over $K$ vanishing at every $P_i$.

If no nonzero such form exists, then $\phi_K$ is injective. We may then select $r(2,D)$ points whose evaluation map is an isomorphism over $K$, retaining the notation $\phi$ for their integral evaluation map. For this selected map, $\det\phi=\bigwedge^{r(2,D)}\phi$ is a nonzero homomorphism between invertible $\O_K$-modules, and we put
\[
 \ell_\fp=\length_{\O_{K,\fp}}
 \coker\bigl(\phi\otimes_{\O_K}\O_{K,\fp}\bigr),\ \fp\in\spm\O_K.
\]
This is the valuation of $\det\phi$ in local lattice bases.

\begin{prop}\label{prop:slope}
Let $Z\subset\PP(\sE_K)(K)$ consist of points of relative height at most $B\geqslant1$. Fix $N_0>1$, and let $\mathcal S$ be a finite set of prime ideals $\fp\in\spm\O_K$ satisfying $\Np\geqslant N_0$. Fix integers $D\geqslant1$ and $h\geqslant0$. For each $\fp\in\mathcal S$, let $T_\fp\subset\PP^2_{\FF_\fp}$ be the zero locus of a nonzero homogeneous form over $\FF_\fp$ of degree at most $D$ (with equation $1$ for the empty set). Suppose that for each $P\in Z$ there are at most $h$ primes $\fp\in\mathcal S$ for which $P_\fp\notin T_\fp$. If
\begin{equation}\label{eq:rank}
 (\#\mathcal S-r(2,D)h)\log N_0>r(2,D)D\log B+\frac{[K:\Q]r(2,D)}{2}\log r(2,D),
\end{equation}
then there exists a nonzero degree $D$ form over $K$ that vanishes on $Z$.
\end{prop}
\begin{proof}
Suppose there is no such form of degree $D$. Then the evaluations on $Z$ span the dual of $E_D\otimes_{\O_K} K$ after trivializing their target lines. Choose $r(2,D)$ independent evaluations, and use $\phi$ and $\ell_\fp$ as defined above. Applying the slope equality of \cite[Proposition~2.1(iii)]{Chen1} to the isomorphism $\phi_K$ gives
\begin{eqnarray*}
 \sum_{\fp\in\spm\O_K}\ell_\fp\log\Np&\leqslant&r(2,D)D\log B-[K:\Q]r(2,D)\,\wmu(\overline E_D)\\
                    &\leqslant&r(2,D)D\log B+\frac{[K:\Q]r(2,D)}{2}\log r(2,D) .  
\end{eqnarray*}
Indeed, the $i$-th target line has Arakelov degree $D[K:\Q]h(P_i)\leqslant D\log B$, and the target is the orthogonal sum of $r(2,D)$ such lines. At each infinite place, each evaluation row has norm at most one because the John norm dominates the supremum norm. Hadamard's inequality therefore bounds the archimedean determinant norms by one. The last inequality follows from \cite[Propositions~2.7 and~2.9, (11) and~(13)]{Chen1} for the trivial rank-three ambient bundle, whose dual tensor powers are orthogonal sums of trivial Hermitian line bundles.

For each of the $r(2,D)$ selected points, at most $h$ primes are exceptional. The union of these exceptional sets therefore has at most $r(2,D)h$ elements. At every remaining $\fp\in\mathcal S$, an equation of $T_\fp$, multiplied by a coordinate monomial of complementary degree, gives a nonzero form of degree $D$ in the kernel of $\phi$ modulo $\fp$. The reduction of $\phi$ is therefore not injective between two vector spaces of dimension $r(2,D)$. Its determinant vanishes modulo $\fp$, so $\ell_\fp\geqslant1$. All other $\ell_\fp$ are nonnegative by integrality, so the left-hand side is at least $(\#\mathcal S-r(2,D)h)\log N_0$, contradicting \eqref{eq:rank}.
\end{proof}

Let $X$, $\delta$ and $\cX$ be as in \S\ref{fundamental setting}. Over each $\O_{K,\fp}$, the closure is defined by a primitive homogeneous equation of degree $\delta$. Its reduction is nonzero, so $\cX_{\FF_\fp}$ is a plane divisor of degree $\delta$. In particular, its local rings at closed points are one-dimensional and Cohen--Macaulay, and $1\leqslant\mu_\xi(\cX_{\FF_\fp})\leqslant\delta$.

For an $\FF_\fp$-rational point $\xi$ of the special fiber $\cX_{\FF_\fp}$ and an integer $r\geqslant1$, we use $Q_\xi(r)$ with the convention of \cite[Notation~15]{Chen2}. It is the sum of the first $r$ terms of the nondecreasing sequence in which $j\geqslant0$ occurs $\dim_{\FF_\fp}(\mathfrak m^j_\xi/\mathfrak m^{j+1}_\xi)$ times, where $\mathfrak m_\xi$ is the maximal ideal of $\O_{\cX_{\FF_\fp},\xi}$. For the general hypersurface estimate see also \cite[\S2.2, (5)--(6)]{Liu2022d}, whose sum at index $r-1$ has $r$ terms.

\begin{lemm}\label{lem:tube}
Let $\delta\geqslant2$, $B\geqslant1$, and let $\xi_i\in\cX_{\FF_{\fp_i}}(\FF_{\fp_i})$ be prescribed at distinct prime ideals $\fp_1,\ldots,\fp_l$, with $l\geqslant1$. If
\begin{equation}\label{eq:tube}
 \sum_{i=1}^l\frac{\log N(\mathfrak p_i)}{\mu_{\xi_i}(\cX_{\FF_{\fp_i}})}
       \geqslant\frac4\delta(\log B+[K:\Q]\log3),
\end{equation}
then all $P\in S(X;B)$ with $P_{\fp_i}=\xi_i$ for every $i$ lie on a curve of degree $\delta-1$ over $K$.
\end{lemm}
\begin{proof}
Put $D=\delta-1$ and $r=r_1(D)$. The rank formula in \S\ref{fundamental setting} gives
\[
r=r_1(D)=r(2,\delta-1)=\binom{\delta+1}{2},\qquad
 \frac{r-1}{2D}=\frac{\delta+2}{4}.
\]

For the Cohen--Macaulay local rings considered above, \cite[Proposition~3.5(ii), (14)]{Chen2} gives \[Q_{\xi_i}(r)\geqslant \frac{r(r-1)}{2\mu_{\xi_i}(\cX_{\FF_{\fp_i}})}.\] We apply \cite[Theorem~3.1, (8)]{Chen2} directly in degree $D$, with all congruence exponents equal to one. It suffices that
\[
 \frac{1}{Dr}\sum_{i=1}^l Q_{\xi_i}(r)\log N(\mathfrak p_i)
 > \log B-\frac{[K:\Q]}{D}\wmu(\overline F_D)+\frac{[K:\Q]\log r}{2D}.
\]
Here we used $[K:\Q]h(P)\leqslant\log B$.

By \cite[\S2.4, (5)]{Chen2} and $r\leqslant3^D$, the right-hand side is at most $\log B+[K:\Q]\log3$. The exact rank and \eqref{eq:tube} give
\begin{align*}
 \frac{1}{Dr}\sum_{i=1}^l Q_{\xi_i}(r)\log N(\fp_i)
 &\geqslant \frac{\delta+2}{4}\sum_{i=1}^l\frac{\log N(\fp_i)}{\mu_{\xi_i}(\cX_{\FF_{\fp_i}})}\\
 &\geqslant \frac{\delta+2}{\delta}(\log B+[K:\Q]\log3)
  >\log B+[K:\Q]\log3 .
\end{align*}
The last inequality is strict even at equality in \eqref{eq:tube}, because $[K:\Q]\log3>0$. The criterion gives the required form.
\end{proof}

\section{Choice of prime ideals and proof of the main theorem}
\subsection{Prime ideals and auxiliary coverings}
Let $X$ be an integral curve in $\mathbb P(\sE_K)$ of degree $\delta\geqslant2$ and $e<B<e^{\delta^2}$. Thus $1<\sqrt{\log B}<\delta$. Put
\begin{equation}\label{eq:parameters}
\begin{aligned}
 &a=128([K:\Q]+1),\qquad N_0=e^{a\sqrt{\log B}},\\
 &M=\lceil(1+\sqrt{\log B})^6\rceil,\qquad r=\lceil\sqrt{\log B}/16\rceil,\\
 &m=\lceil4\delta/\sqrt{\log B}\rceil.
\end{aligned}
\end{equation}

\begin{lemm}\label{lem:primes}
There is a set $\mathcal S$ of $M$ distinct prime ideals of $\OO_K$ satisfying
\begin{equation}\label{eq:norm}
 N_0<\Np\leqslant (2N_0)^{[K:\Q]},\ \forall\;\fp\in\mathcal S.
\end{equation}
\end{lemm}
\begin{proof}
By \cite[Corollary~3, (3.8)]{RS62}, the interval $]N_0,2N_0]$ contains more than $3N_0/(5\log N_0)$ primes in $\mathbb Z$, since $N_0>e^{256}>20.5$. From $e^{a\sqrt{\log B}}\geqslant(a\sqrt{\log B})^8/8!$, $a\geqslant256$ and $\sqrt{\log B}>1$, we get
\[
 \frac{3N_0}{5\log N_0}\geqslant\frac{3a^7(\log B)^{7/2}}{5\cdot8!}
 >128(\log B)^3\geqslant2(1+\sqrt{\log B})^6\geqslant M .
\]
Choose $M$ such primes $p$ and a prime ideal $\fp\in\spm\O_K$ above each. The resulting ideals are distinct, and $\Np=p^{[\FF_\fp:\mathbb F_p]}$ with residue degree between $1$ and $[K:\Q]$, giving \eqref{eq:norm}.
\end{proof}

Fix a set $\mathcal S$ as in Lemma~\ref{lem:primes}. At each $\fp\in\mathcal S$, let $T_\fp$ be the reduced union of the geometric components of the special fiber $\cX_{\FF_\fp}$ occurring with multiplicity at least $m/2$. By Proposition~\ref{lem:mult}, it is defined over $\FF_\fp$, has degree at most $\lfloor2\delta/m\rfloor\leqslant\lfloor \sqrt{\log B}/2\rfloor$, and contains all but at most $\lfloor2\delta^2/m^2\rfloor\leqslant \log B/8$ points of multiplicity at least $m$. For both applications of Lemma~\ref{lem:tube}, we use
\begin{equation}\label{eq:prime-threshold}
 \frac{a\log B}{16}\geqslant4(\log B+[K:\Q]\log3),
\end{equation}
since $(a/64-1)\log B=(2[K:\Q]+1)\log B>[K:\Q]\log3$.

\begin{lemm}\label{lem:small}
The points $P\in S(X;B)$ with $\mu_{P_\fp}(\cX_{\FF_\fp})<m$ at some $\fp\in\mathcal S$ are covered by at most $3M(2N_0)^{2[K:\Q]}$ curves of degree $\delta-1$ in $\mathbb P(\sE_K)$ over $K$. If $\sqrt{\log B}<4$, then these are all the points of $S(X;B)$.
\end{lemm}
\begin{proof}
A residue multiplicity $1\leqslant\mu<m$ satisfies $\mu<4\delta/\sqrt{\log B}$, so $\log\Np/\mu>a\log B/(4\delta)$. By \eqref{eq:prime-threshold}, Lemma~\ref{lem:tube} gives one curve for each residue class of such multiplicity. 

Each prime contributes at most $\#\PP^2(\FF_\fp)=(\Np)^2+\Np+1\leqslant3(2N_0)^{2[K:\Q]}$ classes. 

If $\sqrt{\log B}<4$, we have $m>\delta\geqslant\mu_{P_\fp}(\cX_{\FF_\fp})$ at every prime.
\end{proof}

\begin{lemm}\label{lem:outside}
Suppose $\sqrt{\log B}\geqslant4$. Among the points $P\in S(X;B)$ whose reductions have multiplicity at least $m$ at every selected prime, those with $P_\fp\notin T_\fp$ for at least $r$ primes are covered by at most $(\log B)^{\sqrt{\log B}+1}$ curves of degree $\delta-1$ over $K$.
\end{lemm}
\begin{proof}
For a point in this set, choose any $r$ primes at which $P_\fp\notin T_\fp$ and consider the tuple of its reductions at these primes. The weighted sum in Lemma~\ref{lem:tube} is at least $r\log N_0/\delta\geqslant a\log B/(16\delta)$, since every residue multiplicity is at most $\delta$. Thus one curve suffices for each tuple by \eqref{eq:prime-threshold}.

There are at most $\log B/8$ choices per prime, and $\binom Mr$ choices of the set of primes. A point may occur for more than one choice, which only increases the upper bound. The total number of tuples is at most
\begin{equation}\label{eq:tuples}
 \binom Mr(\log B/8)^r\leqslant(16(\log B)^4)^r\leqslant (\log B)^{5r}\leqslant (\log B)^{\sqrt{\log B}+1}.
\end{equation}
Here $M\leqslant128(\log B)^3$ and $16\leqslant \log B$. For $4\leqslant \sqrt{\log B}\leqslant16$ we have $r=1$ and $5\leqslant\sqrt{\log B}+1$. For $\sqrt{\log B}>16$, use $5r\leqslant5\sqrt{\log B}/16+5<\sqrt{\log B}+1$.
\end{proof}

\begin{lemm}\label{lem:components}
Suppose $\sqrt{\log B}\geqslant4$. The points $P\in S(X;B)$ with $P_\fp\notin T_\fp$ at fewer than $r$ primes in $\mathcal S$ lie on one curve over $K$ of degree $D=\lfloor \sqrt{\log B}/2\rfloor<\delta$.
\end{lemm}
\begin{proof}
Apply Proposition~\ref{prop:slope} in degree $D$, with $s=r(2,D)$ and $h=r-1$. Every point under consideration has at most $h$ exceptions among the $M$ primes. Each $T_\fp$ has an equation of degree at most $D$, including the equation $1$ for the empty divisor. 

Since $s\leqslant(1+\sqrt{\log B})^2$ and $h<\sqrt{\log B}/16$, we have $sh<(1+\sqrt{\log B})^3/16\leqslant M/16$. Moreover $D\leqslant \sqrt{\log B}/2$ and $\log s\leqslant2\sqrt{\log B}$, so
\begin{eqnarray*}
  & &sD\log B+\frac{[K:\Q]s}{2}\log s\\
 &\leqslant& (1+\sqrt{\log B})^2((\sqrt{\log B})^3/2+[K:\Q]\sqrt{\log B})\\
 &\leqslant& ([K:\Q]+1)\sqrt{\log B}(1+\sqrt{\log B})^4\\
 &<&\frac{15}{16}M\log N_0.
\end{eqnarray*}
The last inequality uses $M\geqslant(1+\sqrt{\log B})^6$ and $\log N_0=128([K:\Q]+1)\sqrt{\log B}$. Hence $sD\log B+([K:\Q]s/2)\log s<(M-sh)\log N_0$, as required.
\end{proof}
\subsection{Proof of the main theorem}
We combine the preceding covering lemmas to prove Theorem~\ref{th:main}. 
\begin{proof}[Proof of Theorem~\ref{th:main}]
Lemma~\ref{lem:small} covers every point of small multiplicity at a selected prime. If any points remain then $\sqrt{\log B}\geqslant4$, and Lemma~\ref{lem:outside} or Lemma~\ref{lem:components} applies according as the number of primes for which $P_\fp\notin T_\fp$ is at least $r$ or less than $r$. The total number of curves is at most $3M(2N_0)^{2[K:\Q]}+(\log B)^{\sqrt{\log B}+1}+1$. Since $M\leqslant128(\log B)^3$,
\begin{equation}\label{eq:absorb}
 3M(2N_0)^{2[K:\Q]}+1\leqslant512(\log B)^3\,2^{2[K:\Q]}\exp\!\left(2[K:\Q]a\sqrt{\log B}\right)
                 \leqslant e^{b_K\sqrt{\log B}}.
\end{equation}
For the last inequality, take logarithms and use $\log\sqrt{\log B}\leqslant\sqrt{\log B}$, $\sqrt{\log B}>1$, and $2[K:\Q]a=256[K:\Q]([K:\Q]+1)$. The remaining coefficient satisfies $\log512+6+2[K:\Q]\log2<16+2[K:\Q]$. This gives the constant $b_K$ of Theorem~\ref{th:main} and proves \eqref{eq:bound}.

Each auxiliary form is nonzero and has degree less than $\delta$, so it cannot vanish identically on the integral curve $X$. Its intersection with $X$ is therefore zero-dimensional. Bézout's theorem bounds its degree, hence its number of $K$-points, by $\delta(\delta-1)$.
\end{proof}
\begin{coro}\label{co:diagonal}
Let $K$ be a number field, retain $b_K$ as in Theorem~\ref{th:main}, and let $X$ be an integral plane curve in $\mathbb P(\sE_K)$ of degree $\delta\geqslant1$. For every real $\varepsilon>0$, we have
\[
 \#S(X;\delta)\leqslant
 \bigl(3+e^{2b_K}+e^{b_K^2/(4\varepsilon)}
                     +e^{27/\varepsilon^3}\bigr)\delta^{2+\varepsilon}.
\]
In particular, $\#S(X;\delta)\ll_{[K:\Q],\varepsilon}\delta^{2+\varepsilon}$, uniformly over number fields of fixed degree, with the explicit constant displayed above.
\end{coro}
\begin{proof}
For $\delta\geqslant3$, put $B=\delta$. Completing the square gives $b_K\sqrt{\log\delta}\leqslant\varepsilon\log\delta+b_K^2/(4\varepsilon)$.

The inequality $\log u\leqslant\sqrt u$ holds for $u\geqslant1$, since $(\log u)/\sqrt u$ has maximum $2/e$ at $u=e^2$. Consequently,
\[
 (\sqrt{\log\delta}+1)\log\log\delta\leqslant4(\log\delta)^{3/4}
 \leqslant\varepsilon\log\delta+27/\varepsilon^3,
\]
where the last inequality follows by maximizing $4u^3-\varepsilon u^4$ at $u=3/\varepsilon$ and then taking $u=(\log\delta)^{1/4}$. Theorem~\ref{th:main} now gives the assertion because $e^{\varepsilon\log\delta}=\delta^\varepsilon$.

For $\delta=2$, enlarge the height bound to $B=e^2$. Then $\sqrt{\log B}=\sqrt2<4$, so Lemma~\ref{lem:small} and \eqref{eq:absorb} give at most $e^{b_K\sqrt2}$ lines. By Bézout these contribute at most $2e^{b_K\sqrt2}\leqslant e^{2b_K}2^{2+\varepsilon}$ points.

For $\delta=1$, choose a nonzero coordinate $x_j$ of a point $P$. The product formula gives
\[
 \prod_{v\in M_K}|x_j|_v^{[K_v:\Q_v]}=1.
\]
Each local factor of $H(P)$ is at least the corresponding factor in this product. If another coordinate is nonzero, the inequality is strict at every infinite place, so $H(P)>1$. Hence only the three coordinate points can have height at most one.
\end{proof}
\section*{Declaration}
\subsection*{Conflict of interest statement}
The author declares that there is no potential conflict of interest.
\subsection*{Ethical statement}
The author declares that this manuscript has not been published elsewhere and is not under consideration by another journal. The author will not submit it elsewhere while it is under review. The contribution of this manuscript is the correction of the published proof within its original framework, together with explicit constants uniform over number fields of fixed degree. The author takes responsibility for the mathematical content and the attribution of the results used. The author will inform the editors immediately if any error is found in the manuscript after submission.
\subsection*{Data availability}
No datasets were generated or analysed in this study. 
\backmatter

\bibliography{liu}

\def\cprime{$'$} \def\cprime{$'$}
\providecommand{\bysame}{\leavevmode ---\ }
\providecommand{\og}{``}
\providecommand{\fg}{''}
\providecommand{\smfandname}{\&}
\providecommand{\smfedsname}{\'eds.}
\providecommand{\smfedname}{\'ed.}
\providecommand{\smfmastersthesisname}{M\'emoire}
\providecommand{\smfphdthesisname}{Th\`ese}
\begin{thebibliography}{1}

\bibitem{AM_commutative_algebra}
{\scshape M.~F. Atiyah {\normalfont \smfandname} I.~G. Macdonald} --
  \emph{Introduction to commutative algebra}, Addison-Wesley Publishing Co.,
  Reading, Mass.-London-Don Mills, Ont., 1969.

\bibitem{BCK2025}
{\scshape G.~Binyamini, R.~Cluckers {\normalfont \smfandname} F.~Kato} -- {\og
  Sharp bounds for the number of rational points on algebraic curves and
  dimension growth, over all global fields\fg}, \emph{Proc. Lond. Math. Soc.
  (3)} \textbf{130} (2025), no.~1, p.~Paper No. e70016, 20.

\bibitem{Chen1}
{\scshape H.~Chen} -- {\og Explicit uniform estimation of rational points {I}.
  {E}stimation of heights\fg}, \emph{Journal f\"ur die Reine und Angewandte
  Mathematik} \textbf{668} (2012), p.~59--88.

\bibitem{Chen2}
\bysame , {\og Explicit uniform estimation of rational points {II}.
  {H}ypersurface coverings\fg}, \emph{Journal f\"ur die Reine und Angewandte
  Mathematik} \textbf{668} (2012), p.~89--108.

\bibitem{ARCC02}
{\scshape Colliot-Th{\'e}l{\`e}ne {\normalfont et~al.}} -- {\og Report on the
  {Workshop} `{Rational} and integral points on higher dimensional
  varieties'\fg}, 2002.

\bibitem{Fulton2}
{\scshape W.~Fulton} -- \emph{Algebraic curves. {A}n introduction to algebraic
  geometry}, W. A. Benjamin, Inc., New York-Amsterdam, 1969, Notes written with
  the collaboration of Richard Weiss, Mathematics Lecture Notes Series.

\bibitem{Liu2022d}
{\scshape C.~Liu} -- {\og On the global determinant method\fg}, \emph{Bulletin
  de la Soci\'et\'e{} Math\'ematique de France} \textbf{150} (2022), no.~4,
  p.~699--741.

\bibitem{RS62}
{\scshape J.~B. Rosser {\normalfont \smfandname} L.~Schoenfeld} -- {\og
  Approximate formulas for some functions of prime numbers\fg}, \emph{Illinois
  Journal of Mathematics} \textbf{6} (1962), no.~1, p.~64--94.

\end{thebibliography}
\bibliographystyle{smfplain}

\end{document}